\documentclass[12pt]{amsart}          
\usepackage{amsmath,amssymb}

\theoremstyle{plain}
\newtheorem{xtheorem}{Theorem}
\newtheorem{xlemma}{Lemma}
\newtheorem{xcorollary}{Corollary}
\newtheorem{proposition}{Proposition}

\theoremstyle{remark}
\newtheorem{remark}{Remark}

\DeclareMathOperator{\ch}{ch}
\DeclareMathOperator{\Card}{Card}
\DeclareMathOperator{\height}{ht}

\begin{document}


\title{Graded multiplicities in the exterior algebra}

\thanks{This work was partially supported by Minerva  
grant 8337 and the EC TMR network ``Algebraic Lie Representations''
grant No ERB FMRX-CT97-0100}
\thanks{\copyright{} Academic Press}

\author{Yuri Bazlov}

\address{Department of Mathematics,
The Weizmann Institute of Science,
Rehovot 76100, Israel}

\email{bazlov@wisdom.weizmann.ac.il}

\begin{abstract}This paper 
deals with the graded multiplicities of the ``smallest'' 
irreducible representations of a simple Lie algebra in its exterior algebra.  
An explicit formula for the graded multiplicity 
of the adjoint representation in terms of the Weyl group exponents was 
conjectured by A.~Joseph; a proof of this conjecture, based on the properties 
of Macdonald polynomials, is given in the present paper.
The same method allows to calculate 
the multiplicity of the simple module with highest weight equal to 
the short dominant root.\end{abstract}

\maketitle

\newcommand{\Z}{\mathbb{Z}}     \newcommand{\C}{\mathbb{C}}
\newcommand{\Q}{\mathbb{Q}} 
\newcommand{\RS}{R}             \newcommand{\DRS}{R^\vee}
\newcommand{\RSp}{R^+}          \newcommand{\DRSp}{R^\vee_+}
\newcommand{\Bas}{\Pi}          \newcommand{\DBas}{\dual{\Pi}}           
\newcommand{\rlattice}{Q}	\newcommand{\Drlattice}{\dual{Q}}
\newcommand{\rlatticep}{Q^+}	\newcommand{\Drlatticep}{\dual{Q}_+}
\newcommand{\wlattice}{P}       \newcommand{\Dwlattice}{\dual{P}}       
\newcommand{\wlatticep}{P^+}	\newcommand{\Dwlatticep}{\dual{P}_+}
\newcommand{\spc}{E}             
\newcommand{\dual}[1]{{#1}^\vee}
\newcommand{\W}{W}
\newcommand{\hroot}{\theta}	\newcommand{\Dhroot}{\dual{\hroot}}
\newcommand{\dash}{\nobreakdash-\hspace{0pt}}
\newcommand{\scprod}[2]{\langle{#1},{#2}\rangle_{q,t}}
\newcommand{\cscprod}[2]{({#1},{#2})_{q,t}}
\newcommand{\Pol}{\Q_{q,t}[\wlattice]}
\newcommand{\rank}{r}
\newcommand{\Hecke}{\mathcal{H}_t}
\newcommand{\Neg}{\text{S}}
\newcommand{\hl}{L}
\newcommand{\hs}{S}
\newcommand{\elps}{$\langle\dots\rangle$}
\newcommand{\al}[1]{\alpha^{(#1)}}
\newcommand{\eps}{\epsilon}
\newcommand{\One}[2]{1_{#1}(#2)}
\newcommand{\X}{\frac{t^2-1}{qt^{-2d_{\rank}}-1}}
\newcommand{\g}{\mathfrak{g}}
\newcommand{\h}{\mathfrak{h}}
\newcommand{\q}{\text{GM}}
\newcommand{\exterior}{\Lambda}
\newcommand{\ds}[1]{d^{(s)}_{#1}}


\setcounter{section}{-1}

\section{Introduction}  

\subsection{}

Let $\g$ be a complex simple Lie algebra,
and $V(\lambda)$ a simple $\g$\dash module of highest weight $\lambda$.
The graded multiplicity of $V(\lambda)$ in the exterior algebra 
$\exterior\g$ is a polynomial 
\begin{equation*}
	\q_\lambda(q)=\sum_{n\ge0} [ \exterior^n\g : V(\lambda)] q^n.
\end{equation*}
The module $V(\lambda)$ may appear as an irreducible constituent of
$\exterior\g$, only if $\lambda$ is a dominant weight from the root lattice
of $\g$, lying between $0$ and $2\rho$ in the partial order of weights; here 
$2\rho$ is the sum of the positive roots of $\g$. The polynomials $\q_\lambda$
for such $\lambda$ are not yet known in general, but were calculated for 
some particular values of $\lambda$. For example, if $\lambda$ is close to 
$2\rho$ (exactly: $2\rho-\lambda$ is a combination of simple roots 
with coefficients $0$ and $1$), an explicit formula for 
$\q_\lambda$ is due to Reeder \cite{R}. 
On the other hand, $\q_0$ is the Poincar\'e 
polynomial of the cohomology of the Lie group $G$ corresponding to $\g$,
and is expressed in terms of the exponents 
$d_1,\dots,d_{\rank}$ of the Weyl group of $\g$:
\begin{equation*}
	\q_0(q)=\prod_{i=1}^{\rank}(1+q^{2d_i+1}), 
        \qquad \rank=\text{rank }\g.
\end{equation*}

The above mentioned partial order of weights $\le$ is defined by 
$\lambda\le\mu$, if $\mu-\lambda$ is a sum of positive roots. 
The corresponding Hasse diagram of the dominant part of the root lattice  
looks near zero as 
\begin{equation*}
\begin{cases}
	&0 - \hroot < \hskip-1em - \cdots, \\
	&0 - \hroot_s - \hroot <\hskip-1em - \cdots,
\end{cases}
\end{equation*}
in the case of simply laced and not simply laced Dynkin diagram, respectively;
here $\hroot$ and $\hroot_s$ stand for the highest root and the highest 
short root of $\g$. Thus, $0$, $\hroot$ and $\hroot_s$ are the smallest 
dominant weights, such that the corresponding irreducible 
representations have nontrivial multiplicities in $\exterior\g$. 

In the present paper, we prove the formula for the graded 
multiplicity of the adjoint representation (i.~e. the representation 
$V(\hroot)$): 
\begin{equation*}\label{eq:grmformula}
	\q_{\hroot}(q) = (1+q^{-1}) 
        \bigl(\prod_{i=1}^{\rank-1}(1+q^{2d_i+1})\bigr)
        \sum_{i=1}^{\rank}q^{2d_i},
\tag{$*$}
\end{equation*}
which was conjectured by A.~Joseph in \cite[8.8]{J}. If the Dynkin diagram
is not simply laced, then the representation $V(\hroot_s)$ also appears 
in $\exterior\g$; we calculate its graded multiplicity, and the result is 
\begin{equation*}\label{eq:grmformula1}
	\q_{\hroot_s}(q) = (1+q^{-1}) 
        \bigl(\prod_{i=1}^{\rank-1}(1+q^{2d_i+1})\bigr)
        q^{d_{\rank}+1-2(\rank_s-1)\rank_l}
        \frac{1-q^{4\rank_l\rank_s}}{1-q^{4\rank_l}},
\tag{$**$}
\end{equation*}
where $\rank_s$ is the number of short simple roots 
and $\rank_l$ is the number of 
long simple roots in the root system of $\g$.  

\subsection{}

The method of proving the above formulae was proposed in \cite{J}; 
namely, the graded multiplicity $\q_\lambda$ can be expressed in terms 
of Macdonald polynomials, and it remains to calculate certain value 
of the Macdonald scalar product. In this context, the expression for $\q_0$ 
becomes a direct consequence of Macdonald's constant term formula, 
requiring no cohomological arguments; the computation of the other two 
graded multiplicity polynomials is also carried out in the framework of 
Macdonald theory, involving combinatorial properties of the 
root system and the Weyl group, together with  
arguments from the representation theory of affine Hecke algebras. 

\begin{remark}
The formula for the ungraded multiplicity of the adjoint representation, 
$\q_{\hroot}(1)=2^{\rank}\rank$, which follows from \eqref{eq:grmformula}, 
may be proved using different methods; see \cite[4.2]{R} and references 
therein.
As the author was informed by A.~Joseph, another possible way 
to obtain the graded multiplicity formulae might be to use the 
Clifford algebra techniques from Kostant's paper \cite{Ko}.  
In the case $\g=\mathfrak{sl}(n)$, the graded multiplicities of $V(\lambda)$
where $\lambda$ is a partition of $n$ (in particular, the adjoint 
representation), have been determined combinatorially by Stembridge, see
\cite{S} and \cite[7.4]{R}.
\end{remark}

\subsection{}

The structure of the text is as follows. In section \ref{sect:1} 
we forget about Lie algebras and give some preliminaries 
on root systems and Macdonald polynomials. All the basics on root systems 
can be found in \cite{B}; for the exposition of Macdonald theory and 
Cherednik's proof of Macdonald's conjectures, see 
\cite{C1,Ki,M} and a recent survey \cite{M3}. 
The exponents of the Weyl group are defined purely 
combinatorially in \ref{subsect:exponents}. The next section introduces 
the affine Hecke algebra and its representation in the space of polynomials 
via the Demazure\dash Lusztig operators, which is an important part of 
the current theory of Macdonald polynomials. In this section, we mostly 
adhere to the notation of \cite{M}.
After that we describe the action of certain operator $Y^{\Dhroot}$,
arising from the affine Hecke algebra, on some polynomials. This 
allows to compute particular values of the scalar product $\cscprod{\,}{}$
on the space of polynomials, using unitariness of $Y^{\Dhroot}$ 
with respect to $\cscprod{\,}{}$ and some properties of the scalar product,
proved in section \ref{sect:formulae}. Finally, in section \ref{sect:proof}
we return to the graded multiplicities and give the proof of the formulae 
\eqref{eq:grmformula} and \eqref{eq:grmformula1}, following Joseph's idea 
and using the obtained properties of the scalar product of Macdonald 
polynomials.

It is worth mentioning that we use the Macdonald polynomials depending 
on two parameters $q$, $t$; moreover, all what we need is the case 
$t=q^{-k/2}$ for $k=1$,$2$. In a slightly more general version 
of those polynomials with parameters $q$, $t_l$, $t_s$, the latter two  
are independent when the root system contains both long and 
short roots. We have restricted ourselves to the case $t_l=t_s$, which is 
technically easier (for example, this allows to obtain a rather simple 
scalar product formula valid for both long and short roots, see theorem 
\ref{th:2}) and provides all what we need for the calculation of the 
graded multiplicities.

\subsection{}

The author thanks Professor Anthony Joseph for suggesting the problem 
and giving a number of valuable comments. 


\section{Root system and Macdonald polynomials}
\label{sect:1}

\subsection{Notation}

Let $\RS$ be a reduced irreducible root system spanning 
a finite\dash dimensional vector space $\spc$ of dimension $\rank$ over $\Q$, 
endowed with an inner product $(\cdot\,,\cdot)$.
To each root $\alpha\in\RS$ there corresponds a coroot
$\dual{\alpha}=2\alpha/(\alpha,\alpha)$; the coroots form a dual root 
system $\DRS$.  
We fix a basis $\Bas=\{\alpha_1,\dots,\alpha_\rank\}$ 
of $\RS$ and denote the set of positive roots by $\RSp$, 
the root lattice by $\rlattice$,
the cone spanned by positive roots by $\rlatticep$, the weight lattice by 
$\wlattice$ and the cone of dominant weights by $\wlatticep$. 
We shall also 
use a similar notation $\DBas$, $\DRSp$, $\Drlattice,\dots,\Dwlatticep$ 
for the objects associated with the dual root system.

The orthogonal reflection $s_{\alpha}$ corresponding to a root $\alpha$ acts 
on $\spc$ by 
\begin{equation}\label{eq:s}
	s_\alpha \beta=\beta-(\beta,\dual{\alpha})\alpha, \quad \beta\in\spc.
\end{equation}
The Weyl group $\W$ of $\RS$ is generated by the simple reflections
\begin{equation*}
	s_i=s_{\alpha_i},\quad  i=1,\dots,\rank;
\end{equation*} 
we denote by $\ell(w)$ the length of a reduced decomposition of $w\in\W$ 
with respect to $s_1,\dots,s_{\rank}$.

\subsection{Long and short roots}\label{subsect:2,3}

It is known that all the roots of the same length form 
a $\W$\dash orbit in $\RS$. Furthermore, there are two possible cases: 
\begin{enumerate}
\item
$\RS$ splits into two $\W$\dash orbits, the set of long roots $\RS_l$
and the set of short roots $\RS_s$; there is a number $\Lambda\in\{2,3\}$, 
such that $\frac{(\alpha,\alpha)}{(\beta,\beta)}=\Lambda$ for any 
$\alpha\in\RS_l$, $\beta\in\RS_s$; 
\item 
$\W$ acts transitively on $\RS$, and all the roots are of the same length. 
\end{enumerate}
In the latter case, which corresponds to a simply laced Dynkin diagram
(types $A$, $D$, $E$), 
all the roots are assumed to be long, so $\RS_l=\RS$ and $\RS_s=\emptyset$.
In the former, non\dash simply laced case, $\Lambda=2$ for types $B$, $C$,
$F_4$, $\Lambda=3$ for type $G_2$. 

Since every $\W$\dash orbit meets $\wlatticep$ exactly once, $\RS$ contains 
one or two dominant roots. The long dominant root (the highest root of $\RS$)
will be denoted by $\hroot$; the short dominant root (the highest short root 
of $\RS$), which exists only in non\dash simply laced case, will be denoted by
$\hroot_s$.  

If $\alpha$ and $\beta$ are two roots, then the number 
$(\beta,\dual{\alpha})$ in \eqref{eq:s} may take the following values
\cite[VI, \S1, 3]{B}:

if $\alpha=\pm\beta$, then $(\beta,\dual{\alpha})=\pm2$;

for non\dash proportional $\alpha$, $\beta$, if
$\alpha$ is long or $\beta$ is short, then $(\beta,\dual{\alpha})=0,\pm1$;

if $\alpha$ is short and $\beta$ is long, then 
$(\beta,\dual{\alpha})=0,\pm \Lambda$, where $\Lambda\in\{2,3\}$ is as above. 

\subsection{Weyl group exponents}\label{subsect:exponents}

Recall that $\{\alpha_1,\dots,\alpha_\rank\}$ is the set of simple roots, 
hence a basis of $\spc$. We define a linear function 
$\height:\spc\rightarrow\Q$ by 
\begin{equation*}
	\height \sum_{i=1}^\rank k_i\alpha_i=\sum_{i=1}^{\rank} k_i,
\end{equation*}
so that for a positive root $\alpha$ it gives the height of  
$\alpha$ in the usual sense. Let us introduce two more ``height functions'':
\begin{equation*}
	\height_l \sum_{i=1}^\rank k_i\alpha_i
        =\sum_{\substack{i:\,\alpha_i\text{ is}\\ \text{a long root}}} k_i,
        \qquad
	\height_s \sum_{i=1}^\rank k_i\alpha_i
        =\sum_{\substack{i:\,\alpha_i\text{ is}\\ \text{a short root}}} k_i.
\end{equation*}
Of course, $\height_l+\height_s=\height$; in the simply laced case 
$\height_s=0$.  

For $n$ a positive integer, let, say, $m(n)$ denote the cardinality 
of the set $\{\alpha\in\RSp \mid \height\alpha=n\}$. One can show that 
$m(1)=\rank\ge m(2)\ge\dots$, so those numbers form a partition of 
$\Card\RSp$. The elements of the dual partition, 
$d_1\le d_2\le\dots\le d_\rank$, are called the exponents of the root system 
$\RS$ (or of the Weyl group $\W$). 
As it follows from \cite[V,~\S6,~2]{B}, 
\begin{equation*}
\begin{split} 
	&d_i+d_{\rank+1-i}=d_\rank+1 \quad\text{for }i=1,\dots,\rank, \\
        &d_1=1, \quad d_\rank=\height\hroot.
\end{split}
\end{equation*}
These exponents are determined by the action of $\W$ on the space 
$\spc$ (see \cite[\textit{loc.~cit.}]{B}), so the exponents of the dual root
system are the same, and the partition $\{m(n)\}$ may be defined
via the height function of the dual root root system:
\begin{equation*}
	m(n) = \Card \{\alpha\in\RSp \mid (\rho,\dual{\alpha})=n \}, 
\end{equation*}
where $\rho$ is half the sum of positive roots in $\RS$. 

\subsection{Macdonald polynomials}\label{subsect:macdonald}

Let $\Pol$ be the group algebra of $\wlattice$, generated by 
formal exponentials $e^\lambda$, $\lambda\in\wlattice$ 
over a field $\Q_{q,t}=\Q(q^{1/m},t)$; here $q$ and $t$ are  
two independent variables, and $m$ is an integer such that 
$m\wlattice\subset\rlattice$. We shall refer to the elements of 
$\Pol$ as polynomials. 

The Weyl group $W$ acts on the space of polynomials by 
$w(e^\lambda)=e^{w\lambda}$. Every $W$\dash invariant polynomial is a 
$\Q_{q,t}$\dash linear combination of orbit sums
\begin{equation*}
	m_\lambda=\sum_{\mu\in W\lambda}e^\mu,\qquad \lambda\in\wlatticep.
\end{equation*} 
The Macdonald polynomials are $\W$\dash invariant polynomials 
satisfying the orthogonality and triangularity conditions,
which we now state.

Let us consider the bar involution $f\mapsto \bar{f}$ 
on $\Pol$, which is defined on the exponentials as 
$\overline{e^\lambda}=e^{-\lambda}$ and is extended to $\Pol$ by
$\Q_{q,t}$\dash linearity.
The scalar product related to the Macdonald polynomials may be defined 
as follows: 
\begin{equation}\label{eq:scprod}
      \scprod{f}{g}=\frac{1}{|W|}[f {\bar{g}} \Delta_{q,t}]_0,
\end{equation}
where
\begin{equation*}
      \Delta_{q,t}=\prod_{\alpha\in R}\prod_{i=0}^\infty
\frac{1-q^i e^\alpha}{1-t^{-2}q^i e^\alpha}
\end{equation*}
and $[f]_0$ denotes the constant term, i.~e. the coefficient of $e^0$,  
of a polynomial (or a formal power series) $f$. If the relation 
\begin{equation}\label{eq:k}
	t=q^{-k/2}, \qquad k\text{ is a non-negative integer}, 
\end{equation}
is imposed, $\Delta_{q,t}$ becomes a polynomial 
$\prod_{\alpha\in\RS}(1-e^\alpha)\dots(1-q^{k-1}e^\alpha)$,
so we need no infinite products, cf.\ \cite{M}.
Anyway, the scalar product \eqref{eq:scprod}, which is due to Macdonald 
\cite{M1}, is symmetric and non\dash degenerate. 

We consider a partial order $\le$ on $\wlatticep$: 
\begin{equation}\label{eq:ordering}
	\lambda\le\mu,\text{ if }\mu-\lambda\in\rlatticep.  
\end{equation}
Now we can state the following existence theorem \cite{C1,M,Ki}.
\begin{xtheorem}
\label{th:existence}
There exists a unique family of $\W$\dash invariant polynomials 
$P_\lambda\in\Pol$, $\lambda\in\wlatticep$ such that 

1. $P_\lambda=m_\lambda+\sum_{\mu<\lambda}a_{\lambda\mu}m_\mu$;
\ \ 
2. $\scprod{P_\lambda}{P_\mu}=0$ if $\lambda\ne\mu$.
\end{xtheorem}
These $P_\lambda$'s are called Macdonald polynomials and form a basis 
in the space $\Pol^\W$ of $\W$\dash invariant polynomials.

There is no general formula for the Macdonald polynomials, except for 
some special cases. Namely, if we assume $k=0$ 
in \eqref{eq:k}, then 
$\Delta_{q,t}=1$ and it is easy to see that $P_\lambda=m_\lambda$ satisfy 
the conditions of Theorem \ref{th:existence}. In the case $k=1$
the Macdonald polynomials are given by the Weyl character formula
\begin{equation}\label{eq:Weyl}
	P_\lambda=\chi_\lambda=
{\sum_{w\in\W}(-1)^{\ell(w)}e^{w(\lambda+\rho)-\rho}} \Big/
             {\prod_{\alpha\in\RSp}(1-e^{-\alpha})}
\end{equation}
and are independent of $q$, $t$.


\section{Affine Hecke algebra}

\subsection{Affine root system}

In the notation of the previous section, let us consider the space 
$\widehat{\spc}$ of affine\dash linear functions on $\spc$. The space $\spc$
itself is identified with a subspace of $\widehat{\spc}$ via pairing 
$(\cdot\,,\cdot)$, so $\widehat{\spc}=\spc\oplus\Q\delta$, where $\delta$  
is the constant function $1$ on $\spc$.

We fix the notation for the affine root system 
$\widehat{\RS}=\{\alpha+n\delta \mid \alpha\in\RS, n\in\Z\}
\subset\widehat{\spc}$, the positive affine roots 
$\widehat{\RS}^+=\{\alpha+n\delta \mid \alpha\in\RS,n>0
\text{ or }\alpha\in\RSp,n\ge0\}$,
the zeroth simple affine root $\alpha_0=-\hroot+\delta$, 
and the basis $\widehat{\Bas}=\{\alpha_0\}\cup\Bas$ 
of simple roots in $\widehat{\RS}^+$. 

For an affine root $\widehat{\alpha}$,  
let $s_{\widehat{\alpha}}$ denote the orthogonal reflection of $\spc$
in the affine hyperplane $\{x\mid \widehat{\alpha}(x)=0\}$. 
All $s_{\widehat{\alpha}}$ generate the affine Weyl group $\W^a$, which 
contains $\W$ and is a Coxeter group generated by 
$\{s_0,s_1,\dots,s_{\rank}\}$; here $s_i=s_{\alpha_i}$.

Let $\tau(v)$ denote the translation of $\spc$ by a vector $v\in\spc$. 
Note that 
\begin{equation*}
	\tau(\dual\alpha)=s_{-\alpha+\delta}s_\alpha\in\W^a
	\quad\text{for }\alpha\in\RS,
\end{equation*}
therefore $\W^a$ contains a subgroup of translations $\tau(\Drlattice)$ and 
is a semi\dash direct product  
	$\W^a=\W \ltimes \tau(\Drlattice)$.
The group $\W^a$ is contained in the extended affine Weyl group 
\begin{equation*}
	\widehat{\W}=\W \ltimes \tau(\Dwlattice)
\end{equation*} 
(see \cite[VI, \S2, 3]{B}).
The action of $\widehat{\W}$ on $\spc$, which is, by definition, given by 
$w\tau(\lambda)(x)=w(x+\lambda)$ for $w\in\W$, $\lambda\in\Dwlattice$, 
determines the dual action on $\widehat{\spc}$:
\begin{equation}\label{eq:action}
	w\tau(\lambda)(y+n\delta)=wy+(n-(\lambda,y))\delta,
	\quad w\in\W, \,\lambda\in\Dwlattice, \,y\in\spc, \,n \in\Q.
\end{equation}
This action permutes the affine roots. 

The length function $\ell$ on $\W$ is extended to $\widehat{\W}$ as 
\begin{equation}\label{eq:length}
	\ell(w)=\Card(\widehat{\RS}^+ \cap w^{-1}(-\widehat{\RS}^+)),
	\quad w\in\widehat{\W}.
\end{equation}
The group $\widehat{\W}$ may contain nontrivial elements of zero length, which 
form a finite Abelian group $\Omega$. Each $\omega\in\Omega$ permutes
the simple affine roots $\{\alpha_0,\dots,\alpha_\rank\}$. 

\subsection{Affine Hecke algebra}

The affine Hecke algebra $\Hecke$ of the extended affine Weyl group 
$\widehat{\W}$ with the parameter $t$ 
may be defined in the following way. 
Consider the braid group $B$ of 
$\widehat{\W}$ generated by symbols $T(w)$, $w\in\widehat{\W}$,
subject to the relations 
\begin{equation}\label{eq:relation}
	T(v)T(w)=T(vw),\quad\text{whenever }\ell(v)+\ell(w)=\ell(vw).
\end{equation} 
The algebra $\Hecke$ is the quotient of the group algebra 
$\Q(t)[B]$ modulo the relations 
\begin{equation*}
	(T_i-t)(T_i+t^{-1})=0, \quad i=0,\dots,\rank,
\end{equation*}
where $T_i=T(s_i)$ for a simple reflection $s_i$.

\subsection{Demazure-Lusztig operators}

Now we identify the parameter $t$ in the definition of $\Hecke$ and 
the formal variable $t$ in the space $\Pol$ of polynomials. Besides that, 
we formally put $e^{y+n\delta}=q^{-n}e^y$ for $y+n\delta\in\widehat{\spc}$.
If $w$ is an element of $\widehat{\W}$ decomposed as $v\tau(\lambda)$, 
where $v\in\W$ and $\lambda\in\Dwlattice$, 
then, according to \eqref{eq:action}, 
\begin{equation}\label{eq:action.pol}
   we^\mu=e^{w\mu}=q^{(\lambda,\mu)}e^{v\mu}\quad\text{for }\mu\in\wlattice.
\end{equation} 
(Note that $(\lambda,\mu)$ is always a multiple 
of $1/m$, for $m$ introduced in the definition of $\Q_{q,t}$, see Section 
\ref{sect:1}.) 
This allows to introduce the action of $\Hecke$ on $\Pol$ via the 
Demazure-Lusztig operators 
\begin{equation}\label{eq:DL}
	T_i = t s_i 
            + (t-t^{-1})\frac{1-s_i}{1-e^{\alpha_i}},
         \qquad i=0,\dots,\rank,
\end{equation}
and the rule $\omega e^\mu=e^{\omega\mu}$ for $\omega\in\Omega$
(see \cite{C1} for details). 

\subsection{Operators $Y^\lambda$}

The affine Hecke algebra $\Hecke$ contains a large commutative subalgebra, 
which is important for what follows. Define 
\begin{equation*}
	Y^\lambda=T(\tau(\lambda)),\quad \lambda\in\Dwlatticep.
\end{equation*}
Using \eqref{eq:length} and \eqref{eq:relation}, one can show that 
$Y^{\lambda+\mu}=Y^\lambda Y^\mu$ for dominant $\lambda$, $\mu$. Therefore
the rule $Y^{\nu-\mu}=Y^\nu (Y^\mu)^{-1}$ defines $Y^\lambda$ for arbitrary 
$\lambda\in\Dwlattice$, and $Y^\lambda$ generate a commutative subalgebra 
of $\Hecke$. The algebra $\Hecke$ is generated by $T_1,\dots,T_r$ and 
$Y^\lambda$, $\lambda\in\Dwlattice$.

\subsection{Cherednik's scalar product}\label{subsect:properties}

The scalar product \eqref{eq:scprod} may be replaced by another one,
due to Cherednik \cite{C1}, and the replacement does not affect 
theorem-definition \ref{th:existence} of Macdonald polynomials. Consider 
a new involution $\iota:\Pol\rightarrow\Pol$, which in restriction to 
$\Q_{q,t}$ coincides with the automorphism 
$q\mapsto q^{-1}$, $t\mapsto t^{-1}$, and leaves $e^\lambda$ untouched.
Cherednik's scalar product $\cscprod{\ \,}{\ }$ is defined as follows: 
\begin{equation*} 
 	\cscprod{f}{g}=[f{\bar{g}}^\iota C_{q,t}]_0,
\end{equation*}
where $[\,]_0$, as above, denotes the constant term, and 
\begin{gather*}
 	C_{q,t}=\prod_{\alpha\in R}\prod_{i=0}^\infty
\frac{q^{-(i+\chi(\alpha))/2}e^{\alpha/2}-q^{(i+\chi(\alpha))/2}e^{-\alpha/2}}
{t q^{-(i+\chi(\alpha))/2}e^{\alpha/2}
-t^{-1} q^{(i+\chi(\alpha))/2}e^{-\alpha/2}},
\\
\chi(\alpha)=\begin{cases}
0,\quad \alpha\in\RSp,\\1,\quad-\alpha\in\RSp.
\end{cases}
\end{gather*}
We mention here some properties of this scalar product. First, 
\begin{equation*}
	\cscprod{f}{g}=\cscprod{g}{f}^\iota
	=\cscprod{{\bar f}^\iota}{{\bar g}^\iota}^\iota 
	\quad\text{for }f,g\in\Pol,
\end{equation*}
as it follows from the observation that $\overline{C}_{q,t}=C_{q,t}^\iota$.
Second and important for our setting is that the action 
of $T(w)$, $w\in\widehat{\W}$ is unitary with respect to Cherednik's 
scalar product (see \cite{C1,M,Ki}):
\begin{equation*}
	\cscprod{T(w)f}{g}=\cscprod{f}{T(w)^{-1}g}\quad\text{for }f,g\in\Pol,
\end{equation*}
whence
\begin{equation*}
	\cscprod{Y^\lambda f}{g}=\cscprod{f}{Y^{-\lambda}g}
\quad\text{for all $\lambda\in\Dwlattice$.}
\end{equation*}



\section{Some explicit formulae for the action of the affine Hecke algebra}

\subsection{Action of $Y^{\dual{\hroot}}$}

Our main tool in computing particular values of Cherednik's scalar product  
$\cscprod{\ \,}{\ }$ will be the operator $Y^{\dual{\hroot}}$ on 
$\Pol$, which arises from the action of the affine Hecke algebra, 
as defined in the previous section. 

Recall that $\Lambda\in\{2,3\}$ equals 
$\frac{(\alpha,\alpha)}{(\beta,\beta)}$ for a long root $\alpha$ and a short 
root $\beta$. If $\RS_s=\emptyset$, we assume that the value of $\Lambda$ is 
indeterminate.
The following key Proposition describes the action of $Y^{\dual{\hroot}}$
on certain elements of $\Pol$. It covers the case when the Dynkin diagram 
is simply laced, or $\Lambda=2$; 
the remaining special case $\Lambda=3$, when the root system is 
of type $G_2$, will be considered later. 
\begin{proposition} 
Assume that the root system is not of type $G_2$.
Let $\RSp_s(\hroot)$ denote the set of all positive short roots 
not orthogonal to $\hroot$ (in the simply laced case, this set is
empty). Let us introduce two constants, $\hl=\height_l\hroot$ and 
$\hs=\height_s\hroot$. Then the following holds:
\begin{align}
       &Y^{\Dhroot}e^0=t^{2\hl+\hs}e^0, \label{eq:e0} 
\\      
	\begin{split}
             &Y^{\Dhroot}e^{\hroot}=q^2 t^{-(2\hl+\hs)} e^{\hroot} \\
             &\quad -(t-t^{-1})qt^{-(\hl+\hs)}
              \sum_{\alpha\in\RSp_s(\hroot)}e^\alpha
              -(t-t^{-1})t^{-\hs+1}(qt^{-2\hl}+1)e^0,	
	\end{split}
\notag\tag{\theequation a}\label{eq:ehroot}
\intertext{and, in the non\dash simply laced case,}
    &Y^{\Dhroot}e^{\hroot_s}=qt^{-\hs}e^{\hroot_s}-(t-t^{-1})t^{\hl-\hs+2}e^0.
\notag\tag{\theequation b}\label{eq:ehsroot}
\end{align}
\end{proposition} 
We shall prove the Proposition after we obtain a number of statements 
reflecting some combinatorial structure of the root system and the Weyl group; 
they give rise to properties of Hecke algebra operators and allow to perform 
explicit calculations.

\subsection{Combinatorics}

Now we present some geometrical and combinatorial properties of 
the root system. 

To each element $w\in\widehat{\W}$ there associated a set 
\begin{equation*}
	\Neg(w)=\widehat{\RS}^+ \cap w^{-1}(-\widehat{\RS}^+)
\end{equation*}
of positive affine roots made negative by $w$. By the definition 
\eqref{eq:length} of the length function, $\ell(w)=\Card\Neg(w)$. 
Given a reduced decomposition 
\begin{equation*}
	w=\omega s_{j_p}\dots s_{j_2}s_{j_1}, 
\end{equation*}
where $\omega\in\Omega$ and $s_{j_i}$ are simple reflections, 
$0\le j_i\le\rank$, one has  
\begin{equation}\label{eq:alpha^i}
	\Neg(w)=\{\alpha^{(i)}\mid i=1,\dots,p\}, \quad
	\alpha^{(i)}=s_{j_1}\dots s_{j_{i-1}}\alpha_{j_i}
\end{equation}
(see \cite[VI, \S1, 6]{B},\cite{C1}); 
$\alpha^{(1)},\dots,\alpha^{(p)}$ is the chain of positive roots made 
negative by $w$. The ordering of 
$\alpha^{(i)}$'s depends on the choice of a reduced decomposition of $w$,  
unlike the set $\Neg(w)$ itself. Note that $\Neg(w)\subset\RSp$ for $w\in\W$. 
Some properties of this set\dash valued function $\Neg$ can be found in 
\cite[VI, \S1, 6 and ex.~16]{B}, \cite[section~1]{C1}.

\begin{xlemma}\label{l:formula}
$(i)$ 
Let $w=s_{j_p}\dots s_{j_1}$ be a reduced decomposition of $w\in\W$;  
$s_{j_i}=s_{\alpha_{j_i}}$, $\alpha_{j_i}\in\Bas$. Let 
$\alpha^{(i)}$, $i=1,\dots,p$, be the elements \eqref{eq:alpha^i} of 
$\Neg(w)$. 
Then for $\beta\in\spc$ 
\begin{equation}\label{eq:formula}
	w\beta=\beta-\sum_{i=1}^p (\beta,\dual{\alpha^{(i)}})\alpha_{j_i}.
\end{equation}
$(ii)$ Let $\beta$ be a positive long root. Then the set  
$\Neg(s_\beta)$ consists of $(2\height_l\beta-1)$ long roots and  
$2\Lambda^{-1}\height_s\beta$ short roots.
\end{xlemma}
\begin{proof}
$(i)$ easily follows from formula \eqref{eq:alpha^i} for $\al i$, 
cf.\ \cite[exercise~3.12]{Ka}. 

$(ii)$ Assume $w=s_\beta$ and rewrite \eqref{eq:formula} as 
\begin{equation*}
	-2\beta=-\sum_{i:\,\alpha_{j_i}\text{ is short}}
                 (\beta,\dual{\alpha^{(i)}})\alpha_{j_i}
		-\sum_{i\ne k:\,\alpha_{j_i}\text{ is long}}
		(\beta,\dual{\alpha^{(i)}})\alpha_{j_i}
		-(\beta,\dual{\beta})\alpha_{j_k}
\end{equation*}
(we used that $\beta\in\Neg(s_\beta)$, hence $\beta=\alpha^{(k)}$ for some 
$k$). Note that $\alpha\in\Neg(s_\beta)$ implies $(\beta,\dual{\alpha})>0$; 
therefore, according to \ref{subsect:2,3}, the latter formula reads 
\begin{equation*}
	2\beta=\sum_{i:\,\alpha_{j_i}\text{ is short}}\Lambda\alpha_{j_i}
               +\sum_{i\ne k:\,\alpha_{j_i}\text{ is long}}\alpha_{j_i}
               +2\alpha_{j_k},
\end{equation*}
and it remains to apply $\height_l$ and $\height_s$ to both sides.
\end{proof}

The following two lemmas are of use only in the non\dash simply laced case.
\begin{xlemma}\label{l:>0}
Let $\beta$ be a long root. Assume that a short root $\alpha$ and
an arbitrary root $\gamma$ satisfy $(\alpha,\dual{\beta})>0$, 
$(\gamma,\dual{\beta})>0$. Then $(\gamma,\dual{\alpha})\ge0$; if 
$\gamma$ is short and $\alpha+\gamma\ne\beta$, then $(\gamma,\dual{\alpha})>0$.
\end{xlemma}
\begin{proof}
The hypothesis implies that 
$(\gamma,\dual{\alpha})=(\gamma,s_\beta\dual{\alpha})+
\Lambda(\gamma,\dual{\beta})\ge(\gamma,s_\beta\dual{\alpha})+\Lambda$.
If $\gamma$ is long, the latter is not less than $-\Lambda+\Lambda=0$.
If $\gamma$ is short and $\gamma\ne -s_\beta\alpha=\beta-\alpha$, 
then $(\gamma,s_\beta\dual{\alpha})+\Lambda\ge -1+\Lambda>0$.
\end{proof}
\begin{xlemma}\label{l:ht}
Let $\alpha$ be a short root not orthogonal to $\hroot$, such that 
$\Lambda\alpha-\hroot\in\RSp$. Then 
$\height_l\alpha=(\height_l\hroot+1)/\Lambda$.
\end{xlemma}
\begin{proof}
Such a root $\alpha$ is positive, hence 
$(\alpha,\dual{\hroot})=1$ by \ref{subsect:2,3}; 
$\Lambda\alpha-\hroot=-s_\alpha\hroot$ is a 
positive long root, and we wish to prove that 
$\height_l(\Lambda\alpha-\hroot)=1$. 
By lemma \ref{l:formula}, $(ii)$,
it is enough to show that the set 
$\Neg(s_{\Lambda\alpha-\hroot})$ contains exactly one long root;  
in other words, $\Lambda\alpha-\hroot$ is the only positive long root
made negative by $s_{\Lambda\alpha-\hroot}$.
 
Suppose the set $\Neg(s_{\Lambda\alpha-\hroot})$ contains a long root $\beta$, 
different from $\Lambda\alpha-\hroot$. 
Then $1=(\beta,\dual{(\Lambda\alpha-\hroot)})=
\Lambda(\alpha,\dual{\beta})-
(\hroot,\dual{\beta})$ (note that 
$(\beta_1,\dual\beta_2)=(\beta_2,\dual\beta_1)$, if $\beta_1$, $\beta_2$ are 
of the same length). Since $\beta$ is positive, $(\hroot,\dual{\beta})\ge0$,
therefore the only possible situation is $(\alpha,\dual{\beta})=1$, 
$(\hroot,\dual{\beta})=(\beta,\dual{\hroot})=\Lambda-1$. But the root 
$s_{\Lambda\alpha-\hroot}\beta=\beta+\hroot-\Lambda\alpha$ is negative,
so $0\ge(\beta+\hroot-\Lambda\alpha,\dual{\hroot})=(\Lambda-1)+2-\Lambda=1$
--- a contradiction, which proves the lemma. 
\end{proof}

\subsection{Operators $G_{\alpha}$}  \label{subsect:G_alpha}

Let $w$ be an element of $\widehat{\W}$. To compute the action of 
$T(w)\in\Hecke$ on $\Pol$, one may take a reduced decomposition 
$w=\omega s_{j_p}\dots s_{j_2}s_{j_1}$ 
($\omega\in\Omega$, $s_{j_i}\in\{s_0,\dots,s_{\rank}\}$); 
by virtue of \eqref{eq:relation}, 
\begin{equation*}
	T(w)=\omega T_{j_p}\dots T_{j_2}T_{j_1},	
\end{equation*}
and the action of $T_{j_i}$ is expressed by Demazure\dash Lusztig formula
\eqref{eq:DL}.  
Following \cite{C1}, we introduce the operators $G_\alpha$, 
$\alpha\in\widehat{\RS}$, by 
\begin{equation}\label{eq:G}
	G_\alpha=t + (t-t^{-1})
                 \frac{s_\alpha-1}{1-e^{-\alpha}},
\end{equation}
so $T_i=s_iG_{\alpha_i}$ for $0\le i\le \rank$. 
These operators also possess the property $wG_\alpha w^{-1}=G_{w\alpha}$, 
therefore the latter formula for $T(w)$ can be rewritten as
\begin{equation}\label{eq:decomp.G}
        T(w)=w G_{\alpha^{(p)}}\dots G_{\alpha^{(2)}}G_{\alpha^{(1)}}, 
\end{equation}
where $\alpha^{(i)}=s_{j_1}\dots s_{j_{i-1}}\alpha_{j_i}$ form the chain 
of positive affine roots made negative by $w$, cf.\ \eqref{eq:alpha^i}. 

We shall also use another formula for $G_\alpha$. Let us introduce a notation
\begin{equation*}
	h=t-t^{-1};\qquad
 	\eps(\alpha,\beta)=\begin{cases}
				-1,&\quad(\alpha,\beta)>0, \\
				1,&\quad(\alpha,\beta)\le0,
			   \end{cases}
	\quad \alpha,\beta\in\spc.
\end{equation*}
Then for $\alpha\in\RS$, $\mu\in\wlattice$ one has 
\begin{equation}\label{eq:formula.G}
	G_\alpha e^\mu = t^{\eps} e^\mu + h\eps 
  \sum_{i=1}^{\hskip-5pt |(\mu,\dual{\alpha})|+\frac{\eps-1}{2} \hskip-5pt}
	e^{\mu+i\eps\alpha}, \quad \eps=\eps(\mu,\alpha). 
\end{equation}

\subsection{Decomposition of $Y^{\Dhroot}$}

Since $Y^{\Dhroot}=T(\tau(\Dhroot))$, we should find a reduced 
decomposition of $\tau(\Dhroot)$. First, note that 
\begin{equation*}
	\tau(\Dhroot)=s_{-\hroot+\delta}s_{\hroot}=s_0s_{\hroot}.
\end{equation*}
Let us fix a reduced decomposition 
\begin{equation}\label{eq:reduced.s}
	s_{\hroot}=s_{j_p}\dots s_{j_1}s_{j_0}s_{j_{-1}}\dots s_{j_{-p}}, 
        \quad 1\le j_i\le\rank,
\end{equation}
of the reflection associated with the highest root $\hroot$. The length 
of a reflection has to be odd, say $2p+1$; for further 
convenience, we let the index $i$ run from $-p$ to $p$.
We claim that 
\begin{equation}\label{eq:reduced}
	\tau(\Dhroot)=s_0s_{j_p}\dots s_{j_{-p}}
\end{equation}
is a reduced decomposition of $\tau(\Dhroot)$; this is equivalent to
\begin{equation*}
	\ell(\tau(\Dhroot))=\ell(s_{\hroot})+1.
\end{equation*}
Let us regard $\ell(w)$ as the cardinality of the set $\Neg(w)$. 
It is easy to check that 
\begin{equation*}
\Neg(s_{\hroot})=\{\alpha\in\RSp \mid (\alpha,\Dhroot)>0\},
\quad \Neg(\tau(\Dhroot))=\Neg(s_{\hroot})\cup \{\hroot+\delta\},
\end{equation*}
so $\Card\Neg(\tau(\Dhroot))=\Card\Neg(s_{\hroot})+1$, q.\ e.\ d.  

Let us fix a chain of positive roots made negative by $\tau(\Dhroot)$, 
according to a reduced decomposition \eqref{eq:reduced} and formula 
\eqref{eq:alpha^i}:
\begin{gather*}
	\alpha^{(-p)}=\alpha_{j_{-p}},
   \quad\alpha^{(-p+1)}=s_{j_{-p}}\alpha_{j_{-p+1}},
   \quad\dots,
   \quad\alpha^{(p)}=s_{j_{-p}}\dots s_{j_{p-1}}\alpha_{j_p};\\
\alpha^{(p+1)}=s_{j_{-p}}\dots s_{j_p}\alpha_0=\tau(\Dhroot)^{-1}s_0\alpha_0
=\tau(-\Dhroot)(\hroot-\delta)=\hroot+\delta.
\end{gather*}
Using \eqref{eq:decomp.G}, we obtain
\begin{equation*}
	Y^{\Dhroot}=\tau(\Dhroot)G_{\hroot+\delta}
                    G_{\alpha^{(p)}}\dots G_{\alpha^{(-p)}}.
\end{equation*}
Computing $Y^{\Dhroot}e^\lambda$ (where $\lambda=0$, $\hroot$ or  
$\hroot_s$), we first apply 
$G_{\alpha^{(p)}}\dots G_{\alpha^{(-p)}}$ to $e^\lambda$ and then apply 
$\tau(\Dhroot)G_{\hroot+\delta}$ to the result.

\subsection{Properties of reduced decomposition of $s_{\hroot}$}

We need some properties of reduced decomposition of $s_{\hroot}$ and 
the roots $\alpha^{(-p)},\dots,\alpha^{(p)}$. For our convenience, we choose 
the reduced decomposition \eqref{eq:reduced.s} to be symmetric, i.~e.
\begin{equation*}
	j_{-i}=j_i,\quad i=0,\dots,p.
\end{equation*}
(Such symmetry appears if we represent $s_{\hroot}$ as 
$w^{-1}s_{j_0}w$, where $w$ is the shortest element of $\W$ mapping $\hroot$ 
to a simple long root, say $\alpha_{j_0}$, and then take 
$s_{j_1}\dots s_{j_p}=s_{j_{-1}}\dots s_{j_{-p}}$ to be a reduced 
decomposition of $w$.) Then it is easy to see that 
\begin{equation}\label{eq:rel}
	\alpha^{(-i)}=-s_{\hroot}\alpha^{(i)},\quad i=-p,\dots,p.
\end{equation}
The next lemma deals with the structure of the chain 
$\alpha^{(-p)},\dots,\alpha^{(p)}$, in the simply laced case or case 
$\Lambda=2$. 


\begin{xlemma}\label{l:main}
Assume that the reduced decomposition of $s_{\hroot}$ is chosen to be 
symmetric;
$\Lambda=2$ or the Dynkin diagram of $\RS$ 
is simply laced. Then:

$(a)$~$\alpha^{(-i)}=\hroot-\alpha^{(i)}$ for $i\ne0$; $\alpha^{(0)}=\hroot$.

$(b)$~If $\alpha^{(i)}$ and $\alpha^{(k)}$ are short, then 
$(\alpha^{(i)},\dual{\alpha^{(k)}})=1$ except for $k=\pm i$; 
$(\alpha^{(i)},\dual{\alpha^{(-i)}})=0$. 

$(c)$~If $\alpha^{(i)}$ is short and $\alpha^{(k)}$ is long, $k\ne 0$, then
either
$\genfrac{}{}{0pt}{}{(\alpha^{(i)},\dual{\alpha^{(k)}})=1}{(\alpha^{(i)},
\dual{\alpha^{(-k)}})=0}$, or, vice versa,
$\genfrac{}{}{0pt}{}{(\alpha^{(i)},\dual{\alpha^{(k)}})=0}{(\alpha^{(i)},
\dual{\alpha^{(-k)}})=1}$.

$(d)$~Recall $\hl=\height_l\hroot$. If $\alpha^{(i)}$ is short, then 
$\height_l\alpha^{(i)}=\frac{1}{2}(L+\text{\textit{sgn }}i)$.

$(e)$~Let $\beta_1,\beta_2\dots,\beta_{\hs}$ be the sequence of short roots
obtained from $\al{-p}$, $\dots$, $\al p$ by dropping all the long
roots. This sequence does not depend on the choice of reduced decomposition 
of $s_{\hroot}$; the short root $\beta_m$ is uniquely determined by two
conditions, $(\beta_m,\Dhroot)>0$ and $\height\beta_m=\frac{\hl-1}{2}+m$.

$(f)$~If short roots exist, the highest short root $\hroot_s$ belongs to 
$\Neg(s_{\hroot})$; $\height\hroot_s=\frac{\hl-1}{2}+\hs$.
\end{xlemma}
\begin{proof}
$(a)$ follows immediately from \eqref{eq:rel}. $(b)$: apply lemma \ref{l:>0}
to $\beta=\hroot$, $\alpha=\alpha^{(i)}$ and $\gamma=\alpha^{(k)}$.

$(c)$: by lemma \ref{l:>0}, $(\alpha^{(i)},\dual{\alpha^{(\pm k)}})\ge0$; 
by $(a)$, $(\alpha^{(i)},\dual{\alpha^{(k)}})+
(\alpha^{(i)},\dual{\alpha^{(-k)}})$ equals 
$(\hroot,\dual{\alpha^{(k)}})=1$, so $(c)$ follows.

$(d)$: let $i>0$, $w=s_{j_{i-1}}s_{j_{i-2}}\dots s_{j_{-p}}$. 
Note that $s_{\alpha^{(i)}}\hroot=\hroot-2\alpha^{(i)}$ is a root. 
Obviously,
$\Neg(w)=\{\alpha^{(-p)},\dots,\alpha^{(i-1)}\}$; this set contains 
$\alpha^{(0)}=\hroot$ but does not contain $\alpha^{(i)}$, therefore 
$w(\hroot-2\alpha^{(i)})$ is negative. But $\hroot-2\alpha^{(i)}$ is orthogonal
to $\hroot$, hence cannot lie in $\Neg(w)\subset\Neg(s_{\hroot})$; it means 
that $\hroot-2\alpha^{(i)}$ is itself negative, so lemma \ref{l:ht}
applied to $\alpha^{(i)}$ gives $\height_l \alpha^{(i)}=\frac{1}{2}
(\height_l\hroot+1)$. If $i<0$, 
$\height_l \alpha^{(i)}=\height_l(\hroot-\alpha^{(-i)})=\frac{1}{2}
(\height_l\hroot-1)$. 

$(e)$: first of all, the number of short roots in $\Neg(s_{\hroot})$ is 
equal to $\hs$ by lemma \ref{l:formula}, $(ii)$. As to the enumeration,
we may say that $\beta_m$ is such a short root $\al i$ that 
$\{\al{-p},\dots,\al{i-1}\}$ contains exactly $m-1$ short roots. Now it is 
clearly enough to show that $\height\beta_m=\frac{\hl-1}{2}+m$; by $(d)$, 
this is equivalent to $\height_s\beta_m=m-\frac{1+\text{\it sgn }i}{2}$.
Applying $\height_s$ to both sides of \eqref{eq:formula} with  
$w=s_{j_{i-1}}\dots s_{j_{-p}}$, $\beta=\al i$, we obtain
\begin{equation*}
	1=\height_s\alpha_{j_i}=
        \height_s\al i-\sum_{-p\le k\le i-1:\,\al k\text{ is short}}
	(\al i,\dual{\al k}).
\end{equation*}
By $(b)$, $(\al i,\dual{\al k})=1$ for short $\al k$, except for $k=-i$. 
The index $k=-i$ occurs if and only if $i>0$, therefore the sum on the left
hand side is equal to $m-1-\frac{1+\text{\it sgn }i}{2}$, where 
$m-1=\Card \{\al{-p},\dots,\al{i-1}\}\cap\RS_s$, so $\al i=\beta_m$
and $\height_s\al i=m-\frac{1+\text{\it sgn }i}{2}$.

$(f)$: it is easy to show that two non\dash zero dominant weights cannot be 
orthogonal, hence $(\hroot_s,\Dhroot)>0$ and $\hroot_s\in\Neg(s_{\hroot})$.
Therefore, $\hroot_s\in\{\beta_1,\dots,\beta_{\hs}\}$. Since $\hroot_s$ is the 
highest short root, 
$\height\hroot_s=\max\height\beta_m=\height\beta_{\hs}=\frac{\hl-1}{2}+\hs
$. 
\end{proof}
Recall the function $\eps(\alpha,\beta)=\pm1$ from \ref{subsect:G_alpha}.
Denote
\begin{equation*}
	\eps_{i,j}=\eps(\al i,\dual{\al j}).
\end{equation*}
The following lemma will be used in the computation of 
$Y^{\dual{\hroot}}e^\lambda$:
\begin{xlemma}\label{l:eps}
Suppose $\al i$ is short, $-p\le m\le n\le p$. Then 
\begin{multline*}
	\sum_{k=m}^n \eps_{i,k}
	=n-m+1+2\bigl(-\sum_{k=m}^n (\al i,\dual{\al k})+\One{m,n}{i}\bigr)
	\\ =n-m+1+2
	\bigl(\height(s_{j_n}\dots s_{j_m}-1)s_{j_{m-1}}\dots s_{j_{-p}}\al i
	 +\One{m,n}{i}\bigr),
\end{multline*}
where $\One{m,n}{i}=1$, if $m\le i\le n$, or 
$0$ otherwise.
\end{xlemma}
\begin{proof}
By lemma \ref{l:main}, $(b)$ and $(c)$, 
$(\al i,\dual{\al k})\in\{0,1+\delta_{i,k}\}$, where $\delta_{i,k}$ is the 
Kronecker symbol. Therefore, 
$\eps_{i,k}=1+2\bigl(-(\al i,\dual{\al k})+\delta_{i,k}\bigr)$. Since
$\sum_{k=m}^n \delta_{i,k}=\One{m,n}{i}$, summation over $m\le k\le n$ gives
the first part of the desired formula. To obtain the second part,
apply \eqref{eq:formula}.
\end{proof}

\subsection{Proof of the Proposition}

We are ready to prove the Proposition. 
It is convenient to write 
$Y^{\dual{\hroot}}=\tau(\dual{\hroot})G_{\hroot+\delta}H$, where 
$H=G_{\alpha^{(p)}}G_{\alpha^{(p-1)}}\dots G_{\alpha^{(-p)}}$.

Lemma \ref{l:formula}, $(ii)$, implies that
\begin{equation}\label{eq:p}
	2p+1=\Card \Neg(s_{\hroot}) =2\hl+\hs-1,
\end{equation}
where $\hl=\height_l\hroot$ and $\hs=\height_s\hroot$ as before.

These are particular cases of formula \eqref{eq:formula.G}:  
\begin{gather*}
  G_{\al k}e^{\hroot}=t^{-1}e^{\hroot}-he^{\al{-k}},\,\al k\in\RS_s;\\
  G_{\al k}e^{\hroot}=t^{-1}e^{\hroot}-\delta_{k,0}he^0,\,\al k\in\RS_l;\\
  G_{\al k}e^{\al i}=t^{\eps_{k,i}}e^{\al i}-\delta_{k,i}he^0, 
  \al i\in\RS_s; \quad
  G_{\al k}e^0=te^0
\end{gather*}
(here $\delta_{k,i}$ is the Kronecker symbol). It follows that 
three subspaces of $\Pol$, 
\begin{equation*}
	V_0=\Q_{q,t} e^0 \subset 
	V_1=\sum_{i:\,\al i\in\RS_s}\Q_{q,t} e^{\al i} + V_0 \subset
	V_2=\Q_{q,t} e^{\hroot} + V_1,
\end{equation*}
are $G_{\al k}$\dash invariant for all $\al k\in\Neg(s_{\hroot})$.

All $G_\alpha$ act on $V_0$ by multiplication by $t$, so 
$G_{\hroot+\delta}He^0=t^{2p+2}e^0$. 
Substituting $2p+1=2\hl+\hs-1$ and applying $\tau(\Dhroot)$, which is identity
on $V_0$, we obtain expression \eqref{eq:e0} for $Y^{\Dhroot}e^0$. 

Let us calculate $He^{\hroot_s}\in V_1$. By lemma \ref{l:main}, $(f)$, 
there exists an index $i$ such that $\hroot_s=\al i$. 
Using the above expression for $G_{\al k}e^{\al i}$, we obtain
\begin{equation*}
	He^{\al i}=\prod_{k=-p}^p t^{\eps_{i,k}} \cdot e^{\al i} 
	-t^{p-i}h\prod_{k=-p}^{i-1} t^{\eps_{i,k}} \cdot e^0
\end{equation*}
(the constant term appears from the action of $G_{\al i}$, then each of 
$G_{\al{i+1}}$, $\dots$, $G_{\al p}$ multiplies it by t). 
By lemma \ref{l:eps}, 
$\sum_{k=-p}^p \eps_{i,k}=2p+1+2\bigl(\height(s_{\hroot}-1)\al i+1\bigr)$; 
here $\height(s_{\hroot}-1)\al i=\height(-\hroot)=-\hl-\hs$ and 
$2p+1=2\hl+\hs-1$. 
In the same way, 
\begin{multline}\label{eq:i-1}
	\sum_{k=-p}^{i-1} \eps_{i,k}
	=p+i+2\height(s_{j_{i-1}}  
	\dots s_{-p}-1)\al i
	\\ =p+i+2\height(\alpha_{j_i}-\al i)
	=p+i+2-2\height\al i.
\end{multline} 
Substituting all these in the last expression for $He^{\al i}$, we obtain that
$He^{\al i}$ is equal to 
$t^{-\hs+1}e^{\al i\hskip-3pt}-t^{2\hl+\hs-2\height\al i}he^0$; 
since $\tau(\Dhroot)G_{\hroot+\delta}$ multiplies $e^{\al i}$ by $qt^{-1}$, 
$e^0$ by $t$,  
\begin{equation*} 
	Y^{\Dhroot}e^{\al i}=\tau(\Dhroot)G_{\hroot+\delta}He^{\al i}
	=qt^{-\hs}e^{\al i}-h\cdot t^{2\hl+\hs+1-2\height\al i}e^0.
\end{equation*}
Formula \eqref{eq:ehsroot} follows by replacing $\height\al i=\height\hroot_s$
by $(\hl-1)/2+\hs$ (lemma \ref{l:main}, $(f)$). 

Our next goal is to find $He^{\hroot}$. 
Let us represent $G_{\al k}G_{\al{k-1}}\dots G_{\al{-p}}e^{\hroot}$ $\in V_2$ 
as 
\begin{equation*}
  A_k e^{\hroot}+\sum_{i:\,\al i\text{ is short}}B^{(i)}_ke^{\al i}+C_k e^0,
	\quad A_k,B^{(i)}_k,C_k\in\Q_{q,t}.
\end{equation*}
Examining the above expressions 
for $G_{\al k}\al i$, we conclude that $A_k=t^{-1}A_{k-1}=t^{-p-1-k}$.
The term containing $e^{\al i}$ results from the action of $G_{\al{-i}}$ 
on $A_{-i-1}e^{\hroot}$; thus, $B^{(i)}_k=0$ for $k<-i$, 
$B^{(i)}_{-i}=-hA_{-i-1}$. We therefore have
\begin{equation*}
	B^{(i)}_p=B^{(i)}_{-i} \prod_{k=-i+1}^p t^{\eps_{i,k}}
	=-h \cdot t^{i-p+\sum_{k=-i+1}^p \eps_{i,k}}.
\end{equation*}
In the next subsection we shall show that $\sum_{k=-i+1}^p \eps_{i,k}$ equals
$-i-p-1+\hl$ (see \eqref{eq:debt}), 
so $B^{(i)}_p=-h\cdot t^{-2p-1+\hl}=-h\cdot t^{-\hl-\hs+1}$.
Note that $B^{(i)}_p$ does not depend of $i$. 

Now let us calculate $C_p$. One constant term appears from the action of 
$G_{\al 0}$ on $A_{-1}e^{\hroot}$ and, after that, is $p$ times multiplied 
by $t$, so at the end it gives $t^p(-h\cdot t^{-p})=-h$.
Besides that, if $i>0$ and $\al i$ is short, $G_{\al i}B^{(i)}_{i-1}e^{\al i}$ 
also gives out a constant equal to $-h B^{(i)}_{i-1}$. This constant is 
$p-i$ times multiplied by $t$, and after the $p$th step becomes 
\begin{equation}\label{eq:const}
	t^{p-i}(-h B^{(i)}_{i-1})
        =-h\cdot t^{p-i}\bigl(B^{(i)}_p/\prod_{k=i}^p t^{\eps_{i,k}}\bigr).
\end{equation}
By lemma \ref{l:eps}, $\sum_{k=i}^p \eps_{i,k}$ equals 
\begin{multline*}
	p-i+3+2\height(s_{\hroot}-s_{j_{i-1}}\dots s_{j_{-p}})\al i
	=p-i+3+2\height(\al i-\hroot-\alpha_{j_i})
	\\ =p-i+1-2(\hl+\hs)+2\height\al i,
\end{multline*}
so \eqref{eq:const} is equal to $h^2t^{\hl+\hs-2\height\al i}$.
Summing up:
\begin{multline*}
	C_p=-h+\sum_{\hskip-7pt i>0:\,\al i\text{ is short}\hskip-7pt}
              h^2t^{\hl+\hs-2\height\al i}
	\\ =-h+h^2\sum_{\hskip-5pt m=\hs/2+1\hskip-5pt}^{\hs}
         t^{\hl+\hs-2(\frac{\hl-1}{2}+m)}
	=-h \cdot t^{-\hs};
\end{multline*}
we used that the short roots $\al i$ for $i>0$ are 
$\beta_{\hs/2+1},\dots,\beta_{\hs}$ from lemma \ref{l:main}, $(e)$,
substituted $\height\beta_m=\frac{\hl-1}{2}+m$ and calculated the sum.

It remains to substitute the expressions for $A_p$, $B^{(i)}_p$, $C_p$ to
\begin{equation*}
	He^{\hroot}=A_p e^{\hroot}
	+\sum_{\al i\in\RSp_s(\hroot)} e^{\al i} +C_p e^0
\end{equation*}
and apply $\tau(\Dhroot)G_{\hroot+\delta}$
(note that $G_{\hroot+\delta}e^{\hroot}=t^{-1}e^{\hroot}-hqe^0$). 
The result is formula \eqref{eq:ehroot} for $Y^{\Dhroot}e^{\hroot}$.

\subsection{}

We are left to compute 
$\sum_{k=-i+1}^p \eps_{i,k}=\sum_{k=-p}^{i-1} \eps_{i,-k}$, 
for $i$: $\al i$ is short. Note that if $k=0$ or $\al k$ is a short root 
not orthogonal to $\al i$, then $\eps_{i,-k}=-\eps_{i,k}-2$, since both 
$\eps$'s equal $-1$; otherwise $\eps_{i,-k}=-\eps_{i,k}$.
As it was shown in the proof of lemma \ref{l:main}, $(e)$, there are 
$\height\al i-\frac{\hl-1}{2}-\frac{1+\text{\it sgn }i}{2}$
short roots not orthogonal to $\al i$ among $\al{-p},\dots,\al{i-1}$.
Therefore 
\begin{equation*}
	\sum_{k=-p}^{i-1} \eps_{i,-k}
        =-\sum_{k=-p}^{i-1} \eps_{i,k}
        -2\bigl(\height\al i-\frac{\hl-1}{2}-\frac{1+\text{\it sgn }i}{2}\bigr)
	-2\cdot 1_{-p,i-1}(0).
\end{equation*}
Note that $\frac{1+\text{\it sgn }i}{2}=1_{-p,i-1}(0)$, as $i\ne0$. 
Substituting expression \eqref{eq:i-1} for $\sum_{k=-p}^{i-1} \eps_{i,k}$, we
obtain 
\begin{equation}\label{eq:debt}
	\sum_{k=-i+1}^p \eps_{i,k}=\sum_{k=-p}^{i-1} \eps_{i,-k}=i-p-1+\hl.
\end{equation}
The Proposition has been proved.
\qed

\subsection{The case $G_2$}

Now we compute the action of $Y^{\dual{\hroot}}$ in the case 
when the root system is of type $G_2$. Let $\alpha_1$ (resp.\ $\alpha_2$) 
be the short (resp.\ long) simple root; we have $\hroot=3\alpha_1+2\alpha_2$ 
and $\hroot_s=2\alpha_1+\alpha_2$. Let us denote the short root 
$\alpha_1+\alpha_2$ by $\beta$, and the long root $3\alpha_1+\alpha_2$ 
by $\gamma$. The reflection $s_{\hroot}$ has 
the reduced decomposition $s_2 s_1 s_2 s_1 s_2$, so 
\begin{equation*}
	Y^{\Dhroot}=\tau(\Dhroot)G_{\hroot+\delta}
                    G_{\gamma}G_{\hroot_s}G_{\hroot}G_{\beta}G_{\alpha_2}.
\end{equation*}
The computation gives:
\begin{align}
       &Y^{\Dhroot}e^0=t^6 e^0, \label{eq:e0G2} 
\\      
	\begin{split}
             &Y^{\Dhroot}e^{\hroot}=q^2 t^{-6} e^{\hroot} 
              -(t-t^{-1})qt^{-3}(e^{\hroot_s}+e^\beta)
              -(t-t^{-1})t^{-1}(e^{\alpha_1}+e^{-\alpha_1})
\\             &\qquad -(t-t^{-1})(qt^{-5}+t^{-1})e^0,	
	\end{split}
\notag\tag{\theequation a}\label{eq:ehrootG2}
\\
    &Y^{\Dhroot}e^{\hroot_s}=qt^{-4}e^{\hroot_s}-(t-t^{-1})te^0.
\notag\tag{\theequation b}\label{eq:ehsrootG2}
\end{align}


\section{Scalar product formulae}\label{sect:formulae}

\subsection{Non-symmetry of Cherednik's scalar product}

Our goal is to calculate the value of Cherednik's scalar 
product $\cscprod{e^\beta}{1}$ for any root $\beta$ (here 
$1=e^0=m_0$). The advantage of Cherednik's scalar product is, 
of course, unitariness of the affine Hecke algebra operators 
$T(w)$, which allows to compute scalar product values explicitly; 
but, unlike Macdonald's scalar product $\scprod{\,}{\,}$,   
Cherednik's scalar product is not symmetric, i.~e. 
$\cscprod{wf}{g}$ is not necessarily equal to $\cscprod{f}{g}$, if
$e\ne w\in\W$. 
The next rather general theorem provides some description of this 
non\dash symmetry.

Recall that $\wlattice$ is the weight lattice of the root system $\RS$. 
We say that a subset $A\subset\wlattice$ is convex, if together with two 
weights $\nu$ and $s_\alpha\nu=\nu+k\alpha$ (where $\alpha$ is a simple root 
and $k=-(\nu,\dual{\alpha})>0$), $A$ contains all the weights 
$\nu+\alpha,\dots,\nu+(k-1)\alpha$. Let us call a weight $\nu\in A$ 
maximal in $A$, if $A$ does not contain any weights 
$s_\alpha\nu$, where $\alpha\in\Bas$ and $-(\nu,\dual{\alpha})>0$. 
(For instance, a dominant weight is necessarily maximal.)
\begin{xtheorem}\label{th:convex}
Assume $A$ is a convex subset of $\wlattice$, and there exists a constant $X$
such that the formula
\begin{equation*}
	\cscprod{e^\mu}{1}=t^{-2\height\mu}X
\end{equation*}
holds for any maximal element $\mu$ of $A$. 
Then this formula holds for all $\mu\in A$. 
\end{xtheorem}
\begin{proof}
Recall the partial order $\le$ on $\wlatticep$; now we can compare
two Weyl group orbits in $\wlattice$ by comparing their unique dominant
weights. For $\mu\in A$, let $d(\mu)$ the least integer such that  
there exists a chain $\mu=\mu_{d(\mu)}, \dots,\mu_1,\mu_0$, 
where $\mu_0$ is maximal in $A$ and $\mu_i=s_{\beta_i}\mu_{i+1}$,
$\beta_i$ being a simple root with $-(\mu_{i+1},\dual{\beta}_i)>0$. 
The theorem can be proved by double induction: first in $\W$\dash orbits with
respect to the ordering $\le$, then in $d(\mu)$ inside an orbit.

Let $\mu\in A$. Suppose the formula is already proved for any $\nu\in A$, 
such that $\W\nu<\W\mu$ for the orbits of $\mu$ and $\nu$, and for 
any $\nu\in \W\mu\cap A$, such that $d(\nu)<d(\mu)$.  
If $d(\mu)=0$, $\mu$ is maximal and the formula holds for $\mu$,
otherwise find a simple root 
$\alpha$ such that $k=-(\mu,\dual{\alpha})>0$, $s_{\alpha}\mu\in A$
and $d(s_\alpha\mu)=d(\mu)-1$.   
The operator $T(s_\alpha)=s_\alpha G_\alpha\in\Hecke$ is orthogonal with 
respect to Cherednik's scalar product, so 
\begin{equation*}
t\cscprod{e^{\mu+k\alpha}}{1}
=\cscprod{e^{\mu+k\alpha}}{T(s_\alpha)^{-1}e^0}
=\cscprod{s_\alpha G_\alpha e^{\mu+k\alpha}}{1}.
\end{equation*} 
By formula \eqref{eq:formula.G}, 
\begin{multline*}
	s_\alpha G_\alpha e^{\mu+k\alpha}
	=s_\alpha\bigl(t^{-1}e^{\mu+k\alpha}
        -(t-t^{-1})\sum_{i=1}^{k-1} e^{\mu+i\alpha} \bigr)
	\\=t^{-1}e^{\mu}-(t-t^{-1})\sum_{i=1}^{k-1} e^{\mu+i\alpha}.
\end{multline*}
Now we substitute the latter to the former. 
The formula holds for
$\nu_i=\mu+i\alpha\in A$, $1\le i< k$ (as $\W\nu_i<\W\mu$, which can be seen 
easily), as well as for $\mu+k\alpha$ (as $d(\mu+k\alpha)=d(\mu)-1$). 
So we obtain 
\begin{multline*}
	t\cdot t^{-2k-2\height\mu}X
	=t^{-1}\cscprod{e^\mu}{1}
       -(t-t^{-1})\sum_{i-1}^{k-1}t^{-2i-2\height\mu}X
	\\ =t^{-1}\cscprod{e^\mu}{1}-t^{-1-2\height\mu}X
	+t^{1-2k-2\height\mu}X,
\end{multline*}
therefore $\cscprod{e^\mu}{1}=t^{-2\height\mu}X$.
\end{proof}
\begin{remark}
As it can be seen from the proof, the symmetric polynomial $1$ 
in the theorem can be replaced by arbitrary symmetric 
(i.~e. $\W$\dash invariant) polynomial. 
However, in the present paper we do not use such a generalization.
\end{remark}

\subsection{Certain scalar product values}

Now we are going to calculate explicitly $\cscprod{e^{\hroot}}{1}$ and 
$\cscprod{e^{\hroot_s}}{1}$. The calculation is based on the expressions
of the action of the orthogonal operator $Y^{\Dhroot}$ on $e^0$, 
$e^{\hroot_s}$ and $e^{\hroot}$.

Assume that the Dynkin diagram is not simply laced.
Due to unitariness of $Y^{\Dhroot}$, in all the cases except $G_2$ we have 
the equality 
$\cscprod{Y^{\Dhroot}e^{\hroot_s}}{1}=\cscprod{e^{\hroot_s}}{Y^{-\Dhroot}1}
=t^{2\hl+\hs}\cscprod{e^{\hroot_s}}{1}$ --- note that 
${\overline{Y^{-\Dhroot}1}}^\iota=t^{2\hl+\hs}$ by \eqref{eq:e0}. 
Substituting expression \eqref{eq:ehsroot} for $Y^{\Dhroot}e^{\hroot_s}$, 
we obtain
\begin{equation*}
	q t^{-\hs}\cscprod{e^{\hroot_s}}{1}
      -(t-t^{-1})t^{\hl-\hs+2}\cscprod{1}{1}
	=t^{2\hl+\hs}\cscprod{e^{\hroot_s}}{1},
\end{equation*}
so $\cscprod{e^{\hroot_s}}{1}=\frac{(t-t^{-1})t^{\hl-\hs+2}}{q t^{-\hs}-
t^{2\hl+\hs}}\cscprod{1}{1}$. Using $\height\hroot_s=(\hl-1)/2+\hs$, 
we rewrite this as
\begin{equation}\label{eq:formula_hsroot}
        \cscprod{e^{\hroot_s}}{1} 
        = t^{-2\height\hroot_s}\X\cscprod{1}{1},
\end{equation}
where $d_{\rank}=\height\hroot=\hl+\hs$ is the highest Weyl group exponent, 
see \ref{subsect:exponents}. 
In the $G_2$\dash case this calculation uses expressions \eqref{eq:e0G2} 
for $Y^{\Dhroot}1$ and \eqref{eq:ehsrootG2} for $Y^{\Dhroot}e^{\hroot_s}$, 
so we obtain $q t^{-4}\cscprod{e^{\hroot_s}}{1}-(t-t^{-1})t\cscprod{1}{1}=
t^6 \cscprod{e^{\hroot_s}}{1}$ and come to the same formula 
\eqref{eq:formula_hsroot}, where $\height\hroot_s=3$ and $d_{\rank}=d_2=5$.

Now let us calculate $\cscprod{e^{\hroot}}{1}$. In the same way as for 
$e^{\hroot_s}$, using unitariness of $Y^{\Dhroot}$ and expression 
\eqref{eq:ehroot}, we obtain
\begin{multline}\label{eq:serv_hroot}
	q^2 t^{-(2\hl+\hs)}\cscprod{e^{\hroot}}{1}
	-(t-t^{-1}) q t^{-(\hl+\hs)}\sum_{\alpha\in\RSp_s(\hroot)}
         \cscprod{e^{\alpha}}{1}
 \\	-(t-t^{-1})t^{-\hs+1}(q t^{-2\hl}+1)\cscprod{1}{1}
        =t^{2\hl+\hs}\cscprod{e^{\hroot}}{1}
\end{multline}
in any case except $G_2$. If the short roots exist, the set 
$\RSp_s(\hroot)$ is non\dash empty, so we need to find 
$\cscprod{e^\alpha}{1}$ for $\alpha\in\RSp_s(\hroot)$. Let us use theorem 
\ref{th:convex}. The set $A=\RSp_s$ is convex, and for its unique 
maximal element $\hroot_s$ there holds formula \eqref{eq:formula_hsroot}. 
We conclude that the same formula holds for any $\alpha\in\RSp_s$, 
i.~e. 
\begin{equation}\label{eq:formula_sroot}
	\cscprod{e^{\alpha}}{1} 
        = t^{-2\height\alpha}\X\cscprod{1}{1},
	\qquad \alpha\in\RSp_s.
\end{equation}
We also know that 
$\{\height\alpha\mid \alpha\in\RSp_s(\hroot)\}=
\{\frac{\hl-1}{2}+1,\dots,\frac{\hl-1}{2}+\hs\}$, 
see lemma \ref{l:main}, $(e)$; so \eqref{eq:serv_hroot} reads 
\begin{multline*} 
	q^2 t^{-(2\hl+\hs)}\cscprod{e^{\hroot}}{1}
	-(t-t^{-1}) q t^{-(\hl+\hs)}
      \sum_{m=\frac{\hl-1}{2}+1}^{\frac{\hl-1}{2}+\hs}
         t^{-2m}\X\cscprod{1}{1}
 \\	-(t-t^{-1})t^{-\hs+1}(q t^{-2\hl}+1)\cscprod{1}{1}
        =t^{2\hl+\hs}\cscprod{e^{\hroot}}{1}.
\end{multline*}
Replacing $\sum_m t^{-2m}$ by $t^{-\hl}\frac{1-t^{-2\hs}}{t-t^{-1}}$,
we obtain
\begin{equation*}
	\frac{\cscprod{e^{\hroot}}{1}}{\cscprod{1}{1}}
        =\frac{qt^{-(2\hl+\hs)}(1-t^{-2\hs})\X
          +(t-t^{-1})t^{-\hs+1}
         (q t^{-2\hl}+1)}{q^2 t^{-(2\hl+\hs)} - t^{2\hl+\hs}}.
\end{equation*}
The right\dash hand side simplifies (!) and gives 
\begin{equation}\label{eq:formula_hroot}
	\cscprod{e^{\hroot}}{1} = 
        t^{-2\height\hroot} \X \cscprod{1}{1}
\end{equation}
(here $\height\hroot=d_{\rank}=\hl+\hs$). 

We have proved \eqref{eq:formula_hroot} in all cases except $G_2$. 
For the $G_2$ case, one obtains from \eqref{eq:ehrootG2} that 
\begin{multline*}
    q^2 t^{-6} \cscprod{e^{\hroot}}{1} 
              -(t-t^{-1})qt^{-3}\cscprod{e^{\hroot_s}+e^\beta}{1}
              -(t-t^{-1})t^{-1}\cscprod{e^{\alpha_1}+e^{-\alpha_1}}{1}
\\            -(t-t^{-1})(qt^{-5}+t^{-1})\cscprod{1}{1}	
    =t^6 \cscprod{e^{\hroot}}{1}.
\end{multline*}
Using \eqref{eq:formula_sroot}, we find 
$\cscprod{e^{\hroot_s}}{1}=t^{-6}X$, $\cscprod{e^{\beta}}{1}=t^{-4}X$
and $\cscprod{e^{\alpha_1}}{1}=t^{-2}X$, where 
$X=\frac{t^2-1}{qt^{-10}-1}\cscprod{1}{1}$;
this formula does not work for negative root $-\alpha_1$, 
but it follows from \ref{subsect:properties} that 
$\cscprod{e^{-\alpha_1}}{1}=\cscprod{e^{\alpha_1}}{1}^\iota=qt^{-10}X$. 
Then the computation of $\cscprod{e^{\hroot}}{1}$ leads to the same formula
\eqref{eq:formula_hroot}, which may be easily checked. 

Now we may apply theorem \ref{th:convex} to the convex set $\RSp$. 
Its maximal elements are $\hroot$ and (in the non\dash simply laced case)
$\hroot_s$. Expressions \eqref{eq:formula_hroot} and \eqref{eq:formula_hsroot} 
mean that the formula 
$\cscprod{e^\alpha}{1}=t^{-2\height\alpha}\X\cscprod{1}{1}$ 
holds for $\alpha=\hroot$ and $\alpha=\hroot_s$, so we arrive to
\begin{xtheorem}\label{th:2}
For any positive root $\alpha$
\begin{equation*} 
	\cscprod{e^{\alpha}}{1}=t^{-2\height\alpha}\X \cscprod{1}{1}.
\end{equation*}
\vskip-\baselineskip\vskip-\normallineskip\vskip-\belowdisplayskip\qed
\end{xtheorem}
\noindent
As $\cscprod{e^{-\alpha}}{1}=\cscprod{e^{\alpha}}{1}^\iota$ (see 
\ref{subsect:properties}) and $(\X)^\iota=q t^{-2(d_{\rank}+1)}\X$, 
we obtain also 
\begin{xcorollary}
For any positive root $\alpha$
\begin{equation*} 
	\cscprod{e^{-\alpha}}{1}
        =q t^{2\height\alpha-2(d_{\rank}+1)}\X \cscprod{1}{1}.
\end{equation*}
\vskip-\baselineskip\vskip-\normallineskip\vskip-\belowdisplayskip\qed
\end{xcorollary}


\section{Proof of the graded multiplicity formula}\label{sect:proof}

\subsection{Idea of proof}

The formula for the graded multiplicity of the adjoint representation in its 
exterior algebra 
was conjectured by A.~Joseph in \cite{J}, and a method to prove 
it was also suggested there. Here we present that method, following \cite{J}. 

Let $\g$ be a semi\dash simple Lie algebra over $\C$, and $\h$ be 
its Cartan subalgebra.
If $M$ is a $\g$\dash module and $\mu\in\h^*$, let $M_\mu$ denote the 
$\mu$\dash weight subspace of $M$, and 
$\ch M=\sum_{\mu\in\h^*}\dim M_\mu e^\mu$ the character of $M$.

Furthermore, let 
$M=\oplus_n M_n$ a graded $\g$\dash module. Suppose that $M$ is a direct 
sum of its weight subspaces: $M=\oplus_n \oplus_{\mu\in\h^*}M_{n,\mu}$, 
and that each $M_{n,\mu}$ is finite\dash dimensional. Define the Poincar\'e 
polynomial of $M$ by 
\begin{equation*}
	P_M(q) = \sum_n \sum_{\mu\in\h^*} \dim M_{n,\mu} e^\mu q^n
               =\sum_n \ch M_n q^n.
\end{equation*}

We are interested in the case when $V=\exterior\g$ is the exterior algebra
of the adjoint representation of $\g$, naturally graded by  
$\exterior\g=\oplus_{n\ge0}\exterior^n \g$. It is easy to show that 
\begin{equation*}
	P_{\exterior\g}(q) = (1+q)^{\rank}\prod_{\alpha\in\RS}(1+qe^\alpha),
\end{equation*}
where $\RS$ is the root system of $\g$, and $\rank=\dim\h$ is the rank of 
$\RS$. 

Assume that the algebra $\g$ is simple, so $\RS$ is a reduced irreducible 
root system, and recall the notations $\wlattice$, $\rlattice$ etc.\ from 
Section \ref{sect:1}.  
By graded multiplicity of a simple $\g$\dash module $V(\lambda)$ 
of the highest weight $\lambda\in\wlatticep$ is meant the polynomial 
\begin{equation*}
	\q_\lambda(q)=\sum_{n\ge0} [ \exterior^n\g : V(\lambda)] q^n.
\end{equation*}
This may be expressed in terms of Macdonald theory as follows. 
For $k=0,1,\ldots$, denote by $\langle f,g \rangle_k$ 
the scalar product $[f \bar g \Delta_k]_0$ on $\Pol$, 
where 
\begin{equation*}
	\Delta_k=\prod_{\alpha\in\RS}\prod_{i=0}^{k-1}(1-q^i e^\alpha)
\end{equation*}
is the specialization $t=q^{-k/2}$ of Macdonald's $\Delta_{q,t}$ 
(see \ref{subsect:macdonald}). 
Note that $P_{\exterior\g}(q)$ is equal to $\sum_{\lambda\in\wlatticep}
\q_\lambda(q) \ch V(\lambda)$; the characters 
$\ch V(\lambda)=\chi_\lambda$ are given by the Weyl character formula 
\eqref{eq:Weyl} and form an orthonormal basis of $\Pol$ with respect 
to $\langle\ \,,\ \rangle_1$. Therefore
\begin{equation*}
	\q_\lambda(q) =  \langle P_{\exterior\g}(q), \chi_\lambda \rangle_1.
\end{equation*}
Since $P_{\exterior\g}(-q)=(1-q)^{\rank} \Delta_2/\Delta_1$, 
we have
\begin{equation*}
	\q_\lambda(-q) = (1-q)^{\rank} \langle 1,\chi_\lambda \rangle_2.
\end{equation*}
It is clear enough (e.~g., from this formula) that 
$\q_\lambda(q)=0$, if $\lambda\not\in\rlattice$. So the problem is to find
$\langle 1,\chi_\lambda \rangle_2$ for $\lambda\in\wlatticep\cap\rlattice$. 
In what follows, we calculate this for the smallest dominant elements
of the root lattice, namely for $\lambda=0,\hroot_s,\hroot$; the graded
multiplicity of the adjoint representation in its exterior algebra is 
$\q_{\hroot}$.
 
(On the other hand, the formula for $\lambda$ close to 
$2\rho=\sum_{\alpha\in\RSp}\alpha$ was found in \cite{R}, see Introduction; 
if $\lambda > 2\rho$, $\q_\lambda$ is zero.)
              
\subsection{Calculation}

Let us start with the multiplicity of the trivial representation. 
We have 
$\q_0(-q)=(1-q)^{\rank}\langle 1,1 \rangle_2$. The formula 
for $\langle 1,1 \rangle_2$ is one of (now proved) Macdonald's constant term 
conjectures; from \cite{M}
\begin{equation*}
	\langle 1,1 \rangle_2 = \prod_{\alpha\in\RSp} 
        \frac{1-q^{2(\rho,\dual{\alpha})+1}}{1-q^{2(\rho,\dual{\alpha})-1}}.
\end{equation*}
We rewrite this formula as 
$(1-q)^{-\rank}\prod_{n\ge0} (1-q^{2n+1})^{m(n)-m(n+1)}$, 
where $m(n)$ is the cardinality of 
$\{\alpha\in\RSp\mid(\rho,\dual{\alpha})=n\}$. Then, since 
$d_1,\dots,d_{\rank}$ is the partition dual to $m(1)\ge m(2)\ge\dots$
(recall \ref{subsect:exponents}), we have
\begin{equation*}
	\q_0(-q)=\prod_{i=0}^{\rank} (1-q^{2d_i+1}).
\end{equation*}

Now we proceed with $\q_{\hroot}$.  To find $\langle 1,\chi_{\hroot}\rangle_2$,
note first that the character $\chi_{\hroot}$ 
of the adjoint representation 
is $\rank+\sum_{\alpha\in\RS} e^\alpha$, so
\begin{equation*}
   \frac{\q_{\hroot}(-q)}{\q_0(-q)}=
   \frac{\langle 1,\chi_{\hroot}\rangle_2}{\langle 1,1 \rangle_2} 
        = \rank + 
\frac{\langle 1,\sum_{\alpha\in\RS} e^\alpha \rangle_2}{\langle 1,1 \rangle_2}.
\end{equation*}
By \cite{M}, for any symmetric polynomials $f,g$, 
$\langle f,g \rangle_2=\langle g,f \rangle_2$ equals 
$c_2(g,f^\iota)_2$, where $(\,,)_2$ is the value of Cherednik's scalar product
$(\,,)_{q,t}$ under the relation $t=q^{-1}$, 
and $c_2$ is a constant independent
of $f$, $g$. It is thus enough to find 
$\cscprod{\sum_{\alpha\in\RS}e^\alpha}{1}/\cscprod{1}{1}$. 
By Theorem \ref{th:2} and Corollary, 
\begin{equation}\label{eq:fract}
    \frac{\cscprod{\sum_{\alpha\in\RS}e^\alpha}{1}}{\cscprod{1}{1}}
    =\X \sum_{\alpha\in\RSp}
    (t^{-2\height\alpha}+qt^{2\height\alpha-2(d_{\rank}+1)}). 
\end{equation}
We calculate this with the aid of the following easy\dash to\dash prove lemma,
which follows directly from the definition of dual partition:
\begin{xlemma}\label{l:d}
\begin{equation*}
	\sum_{\alpha\in\RSp}t^{k\height\alpha+l}
	=\sum_{i=1}^{\rank} \frac{t^l}{t^{-k}-1} (1-t^{k d_i}),
\end{equation*}
where $d_1,\dots,d_{\rank}$ are the Weyl group exponents.
\qed
\end{xlemma}
Using the lemma and the relation $d_{\rank}-d_i=d_{\rank+1-i}-1$, 
we rewrite \eqref{eq:fract} as
\begin{multline}\label{eq:scprod_hroot}
   \X\sum_{i=1}^{\rank} \bigr( \frac{1-t^{-2d_i}}{t^2-1}
       +\frac{qt^{-2(d_{\rank}+1)}(1-t^{2d_i})}{t^{-2}-1} \bigr)
    \\ = \sum_{i=1}^{\rank}
      \frac{1-t^{-2d_i} - qt^{-2d_{\rank}}(1-t^{2d_i})}{qt^{-2d_{\rank}}-1}
    = -\rank+\sum_i 
  \frac{q t^{-2(d_{\rank}-d_i)} - t^{-2d_i}}{q t^{-2d_{\rank}}-1}
  \\= -\rank + \frac{q t^2 -1}{q t^{-2d_{\rank}}-1}\sum_i t^{-2d_i}.
\end{multline}
Substituting $t=q^{-1}$, we come to
\begin{equation*}
	\q_{\hroot}(-q) 
       = \q_0(-q)\frac{q^{-1}-1}{q^{2d_{\rank}+1}-1}\sum_{i=1}^{\rank}q^{2d_i},
\end{equation*}
which, together with the expression for $\q_0$, gives the graded multiplicity
formula \eqref{eq:grmformula}.

\subsection{Graded multiplicity of $V(\hroot_s)$}

Finally, we calculate the graded multiplicity polynomial 
$\q_{\hroot_s}(q)$ in the non\dash simply
laced case. Let $\rank_s$ be the number of short simple roots, and 
$\rank_l=\rank-\rank_s$ the number of long simple roots.
For $n\ge0$, let $m_s(n)$ denote the number of 
short positive roots of height $n$. One can show that 
$m_s(1)\ge m_s(2)\ge\dots$; 
let $\ds 1$,$\ds 2,\dots\ds{\rank_s}$ be the partition dual to $\{m_s(n)\}$.
These $\ds i$'s look like some analogue of the exponents 
$d_1,\dots,d_{\rank}$, but they always form an arithmetic progression  
\begin{equation*}
	\ds i = (d_{\rank}+1)/2+(2i-1-\rank_s)\rank_l,
        \qquad i=1,\dots,\rank_s
\end{equation*}
(the simplest way to check this is may be the direct verification; 
note that $\rank_s$ may exceed $2$ in $C_n$ case only). 
In particular, $d_{\rank}-d_i=d_{\rank_s+1-i}-1$. The latter relation 
allows to obtain a formula quite similar to \eqref{eq:scprod_hroot}, 
using theorem \ref{th:2}, corollary and lemma \ref{l:d} applied to $\RSp_s$: 
\begin{equation*}
	\frac{\cscprod{\sum_{\alpha\in\RS_s}e^\alpha}{1}}{\cscprod{1}{1}}
       = -\rank_s + \frac{q t^2 -1}{q t^{-2d_{\rank}}-1}\sum_i t^{-2\ds i}.
\end{equation*}
(The proof is the same as of \eqref{eq:scprod_hroot}.)
Using the expression for $\ds i$, we may write this formula as 
\begin{equation}\label{eq:scprod_hsroot}
	\frac{\cscprod{\sum_{\alpha\in\RS_s}e^\alpha}{1}}{\cscprod{1}{1}}
       = -\rank_s + \frac{q t^2 -1}{q t^{-2d_{\rank}}-1}
         t^{-d_{\rank}-1+2(\rank_s-1)\rank_l}
            \frac{1-t^{-4\rank_l\rank_s}}{1-t^{-4\rank_l}}.
\end{equation}

Now $\ch V(\hroot_s)=\chi_{\hroot_s}$ is the Macdonald polynomial 
$P_{\hroot_s}$ subject to the relation $t=q^{-1/2}$. Note that 
$P_{\hroot_s}=m_{\hroot_s}-\frac{\cscprod{m_{\hroot}}{1}}{\cscprod{1}{1}}e^0$,
where $m_{\hroot_s}=\sum_{\alpha\in\RS_s}$ is the orbit sum of $\hroot_s$.
If we assume $t=q^{-1/2}$ in \eqref{eq:scprod_hsroot}, the second term in the
right\dash hand side vanishes; we get 
$
	\chi_{\hroot_s} = m_{\hroot_s} + \rank_s
$, 
so 
\begin{equation*}
	\frac{\q_{\hroot_s}(-q)}{\q_0(-q)}
        = \frac{(\chi_{\hroot_s},1)_2}{(1,1)_2}
	= \frac{q^{-1}-1}{q^{2d_{\rank}+1}-1}
          q^{d_{\rank}+1-2(\rank_s-1)\rank_l}
         \frac{1-q^{4\rank_l\rank_s}}{1-q^{4\rank_l}},
\end{equation*}
which gives \eqref{eq:grmformula1}. To find 
$\frac{(\chi_{\hroot_s},1)_2}{(1,1)_2}=\rank_s+
\frac{(m_{\hroot_s},1)_2}{(1,1)_2}$, we assumed $t=q^{-1}$ in 
\eqref{eq:scprod_hsroot}.


\end{document}